\documentclass[12pt,a4paper]{article}
\usepackage[latin2]{inputenc}
\usepackage{amsmath}
\usepackage{amsfonts}
\usepackage{amssymb}
\usepackage{graphicx}
\usepackage[left=2.00cm, right=2.00cm, top=2.50cm, bottom=2.50cm]{geometry}
\usepackage[T1]{fontenc}

\def\Q{{\mathbb Q}}
\def\Z{{\mathbb Z}}

\newtheorem{lemma}{Lemma}
\newtheorem{theorem}[lemma]{Theorem}
\newtheorem{conjecture}[lemma]{Conjecture}

\newtheorem{proposition}[lemma]{Proposition}

\title{
On the monogenity of totally complex pure octic fields
}
\author{
Istv\'{a}n Ga\'{a}l \\
{\small University of Debrecen, Mathematical Institute} \\
{\small H--4002 Debrecen Pf.400., Hungary,} \\
{\small e--mail: gaal.istvan@unideb.hu},
}

\begin{document}
\baselineskip=17pt

\maketitle
\thispagestyle{empty}

\renewcommand{\thefootnote}{\arabic{footnote}}
\setcounter{footnote}{0}

\vspace{0.5cm}


\begin{abstract}
Let $0,1\ne m\in\Z$ and $\alpha=\sqrt[8]{m}$. 
According to the results of I. Ga\'al and L. El Fadil \cite{smalldeg}, \cite{octic},
$\alpha$ generates a power integral basis in $K=\Q(\alpha)$, if and only if
$m$ is square-free and $m\not\equiv 1\;(\bmod\; 4)$.
In the present paper we consider totally complex pure octic fields, that is the case $m<0$,
with $m$ satisfiying the above property. In this case $(1,\alpha,\alpha^2,\ldots,\alpha^7)$
is an integral basis. Our purpose is to investigate whether $K$ admits any other 
generators of power integral bases, inequivalent to $\alpha$.
We present an efficient method to calculate generators of power integral bases
in this type of fields with coefficients $<10^{200}$ in the above integral basis.
We report on the results of our calculation for this type of fields with
$0>m>-5000$, which yields 2024 fields.
\end{abstract}


\section{Introduction}

Monogenity properties and the existence of a power integral bases in number fields
is a classical area of algebraic number theory, cf. \cite{nark}, \cite{book}.
The number field $K$, or its ring of integers $\Z_K$ is called {\it monogenic},
if $\Z_K$ is mono-generated, that is $\Z_K=\Z[\alpha]$ with some algebraic integer $\alpha$, 
or equivalently,
if $K$ has a {\it power integral basis} of type $(1,\alpha,\alpha^2,\ldots,\alpha^{n-1})$.

During the last couple of years monogenity properties of several types of number fields
were investigated by several authors. One of the most frequently studied type of number fields are 
the {\it pure fields}, generated by a root $\alpha$ of an irreducible polynomial
$x^n-m$.
In addition to classical methods (see e. g. I. Ga\'al and L. Remete \cite{gr2017}) 
it has turned out, that the Newton polygon method
can be efficiently applied to consider monogenity in pure fields, 
see 
J. Guardia, J. Montes and E. Nart \cite{GMN},
O. Ore \cite{O},
I. Ga\'al and L. El Fadil \cite{smalldeg}.

In a recent paper \cite{ig6} we studied totally complex pure sextic fields and found that
if $\sqrt[6]{m}$ generates a power integral basis in $\Q(\sqrt[6]{m})$, then up to equivalence
(sign translation by elements of $\Z$) this is the only generator of power integral bases, among elements with coordinates in absolute value $<10^{100}$
in the integral basis.

In the present paper we investigate a similar phenomena in totally complex pure octic fields.

In \cite{g42} we gave an algorithm to determine power integral bases in any octic fields with a quadratic subfiled.
The algorithm provides all possible generators, but the calculation is tedious: 
it requires the resolution of a cubic and some corresponding quartic relative Thue equations over the quadratic subfield,
as well as the resolution of a unit equation in two variables in a number field of degree 12.

On the other hand we have the experience that such relative Thue equations have only small solutions,
mostly a few digits in the integral bases. In order to be able to deal with a large number of octic fields 
we apply here an efficient algorithm \cite{small} providing the solutions in absolute values up to $10^{200}$
concerning the coefficients of the variables in an integer basis of the ground field.
Having these solutions we can be very sure that these are indeed all solutions.

For the completeness we add, that monogenity of several other types of octic fields were also
formerly successfully considered, even infinite parametric families of octic fields,
 see \cite{gsz13}, \cite{grsz14}, \cite{grsz16}
\cite{gr2017b}, \cite{gr2018a}.  In these cases the octic fields are mostly
 composites of a quadratic and a quartic field, which makes the resolution much easier.

In the present paper we construct an efficient algorithm applying all special properties 
of pure octic fields. We shall compose the methods of 
I. Ga\'al and M. Pohst \cite{g42} for the calculation of relative power integral bases in quartic relative
extensions and I. Ga\'al \cite{small} for calculating "small" solutions of relative Thue equations.

\section{Totally complex pure octic fields}

Let $0\ne m\in\Z$ be a rational integer, $\alpha=\sqrt[8]{m}$, $K=\Q(\alpha)$
with ring of integers $\Z_K$.
L. El Fadil and I. Ga\'al \cite{smalldeg}, \cite{octic} showed

\begin{proposition}
$\alpha$ generates a power integral basis in $K$, if and only if
\begin{center}
(*) \hspace{3cm} $m\ne 1$ is square-free, $m\not\equiv 1\;(\bmod\; 4)$.
\end{center}
\end{proposition}

In the present paper we develop a fast algorithm for determining generators 
of power integral bases in pure octic with $0>m>-5000$ with property (*).
There are 2024 fields of this type. 

In all these number fields $(1,\alpha,\alpha^2,\alpha^3,\alpha^4,\alpha^5,\alpha^6,\alpha^7)$ is an integral basis,
hence our purpose is to decide if there exist other non-equivalent generators of power integral bases.
($\alpha,\gamma\in\Z_K$ are called equivalent, if $\gamma=a\pm\alpha$. Obviously, $\alpha$ generates a power 
integral basis in $K$ if and only if the equivalent $\gamma$ has this property.)

The above integral basis is 
the same as $(1,\alpha,\alpha^2,\alpha^3,\sqrt{m},\alpha\sqrt{m},\alpha^2\sqrt{m},\alpha^3\sqrt{m})$, hence
we represent any $\gamma\in\Z_K$ in the form
\begin{equation}
\gamma=x_0+x_1\alpha+x_2\alpha^2+x_3\alpha^3+y_0\sqrt{m}+y_1\alpha\sqrt{m}
+y_2\alpha^2\sqrt{m}+y_3\alpha^3\sqrt{m}=W+\alpha X+\alpha^2 Y+\alpha^3 Z,
\label{gg}
\end{equation}
with $x_0,x_1,x_2,x_3,y_0,y_1,y_2,y_3\in\Z$,
\[ 
W=x_0+\sqrt{m}\ y_0,X=x_1+\sqrt{m}\ y_1,Y=x_2+\sqrt{m}\ y_2,Z=x_3+\sqrt{m}\ y_3
\]
algebraic integers in $M=\Q(\sqrt{m})$
having integral basis $(1,\sqrt{m})$ under the condition (*).

We shall determine all generators of power integral bases in $K$ with coordinates 
\begin{equation}
\max(|y_0|,|x_1|,|y_1|,|x_2|,|y_2|,|x_3|,|y_3|)<10^{200}.
\label{100}
\end{equation}

Surprisingly, our calculations show:
\begin{theorem}
In the 2024 fields $K=\Q(\sqrt[8]{m})$ with $0>m>-5000$ satisfying (*), 
up to equivalence $\sqrt[8]{m}$ is
the only generator of power integral bases with (\ref{100}), except for $m=-1$,
when we additionally have $i\sqrt[8]{m},(\sqrt[8]{m})^3,i(\sqrt[8]{m})^3$.
\end{theorem}

Our experience shows that considering a certain type of number fields, the monogenic 
fields are more frequent with smaller absolute discriminants, and in those monogenic fields
generators of power integral bases usually have small coefficients in the integral basis.
Therefore, based on our calculation we have the following

\begin{conjecture}
If $m<0$ satisfies (*) then up to equivalence $\sqrt[8]{m}$ is
the only generator of power integral bases in $K=\Q(\sqrt[8]{m})$, except for $m=-1$,
when we additionally have $i\sqrt[8]{m},(\sqrt[8]{m})^3,i(\sqrt[8]{m})^3$.
\end{conjecture}

\section{Generators of relative and absolute power integral bases}

If $\gamma$ is a generator of a power integral bases of $K$ (in the absolute sense, over $\Q$),
then it is also a generator of a relative power integral basis of $K$ over $M$ (cf. \cite{book}).
Therefore we shall first determine all possible 
\[
\delta=X\alpha+Y\alpha^2+Z\alpha^3
\]
such that $(1,\delta,\delta^2,\delta^3)$ is a relative  integral basis of $K$ over $M$.
These $\delta$ are determined up to $\Z_M$-equivalence, that is up to multiplication by a unit in $M$
and translation by an element of $\Z_M$.
Hence we are looking for generators of (absolute) power integral bases of $K$ in the form
\[
\gamma=W+\varepsilon\cdot \delta
\]
with $W\in\Z_M$, $\varepsilon$ a unit in $M$. We have $W=x_0+\sqrt{m}\, y_0$ and $\gamma$ is also
only determined up to $\Z$-equivalence (up to sign and translation by elements of $\Z$),
therefore in addition to $x_1,y_1,x_2,y_2,x_3,y_3$ in $X,Y,Z$ we additionally only have to determine
$\varepsilon$ and $y_0$. But $M$ has only $\pm 1$ as units ($\pm 1,\pm i$ for $m=-1$), hence
$\varepsilon$ does not make any difficulty and $y_0$ can be determined form the condition that the index of 
of $\gamma$ is $I(\gamma)=1$.

Finally, if $x_1,y_1,x_2,y_2,x_3,y_3$ in $X,Y,Z$ satsify (\ref{100}), then this is equivalently valid
for the coordinates of $\varepsilon X,\varepsilon Y,\varepsilon Z$. In the last steps we
shall see that the corresponding $y_0$ also remains small.

\section{Constructing a quartic relative Thue equation over $M$}

First we formulate a special case of the main result of \cite{g42} (see also \cite{book}) what we shall use in the sequel.
As above up to $\Z_M$-equivalence any generator $\delta\in\Z_K$ of a power integral basis 
of $K$ over $M$ can be written in the form
\begin{equation}
\delta=X\alpha+Y\alpha^2+Z\alpha^3,
\label{teta}
\end{equation}
with $X,Y,Z\in\Z_M$.
Let 
\[
f(x)=x^4+a_1x^3+a_2x^2+a_3x+a_4\in\Z_M[x]
\]
 be the relative minimal polynomial of $\alpha$ over $M$ and let
\[
F(U,V)=U^3-a_2U^2V+(a_1a_3-4a_4)UV^2+(4a_2a_4-a_3^2-a_1^2a_4)V^3
\]
a binary cubic form over $\Z_M$ and
\begin{eqnarray*}
Q_1(X,Y,Z)&=&X^2-XYa_1+Y^2a_2+XZ(a_1^2-2a_2)+YZ(a_3-a_1a_2)+Z^2(-a_1a_3+a_2^2+a_4)\\
Q_2(X,Y,Z)&=&Y^2-XZ-a_1YZ+Z^2a_2
\end{eqnarray*}
ternary quadratic forms over $\Z_M$.

\begin{lemma}(\cite{g42})
If $\delta$ of (\ref{teta}) generates a relative power integral basis
of $K$ over $M$, then 
then there is a solution $(U,V)\in \Z_M$ of
\begin{equation}
N_{M/\Q}(F(U,V))=\pm 1,
\label{res}
\end{equation}
such that
\begin{eqnarray}
U&=&Q_1(X,Y,Z), \nonumber \\
V&=&Q_2(X,Y,Z).   \label{q12}
\end{eqnarray}
\label{l1}
\end{lemma}  

In our case the minimal polynomial of $\alpha=\sqrt[8]{m}$ over $M=\Q(\sqrt{m})$ is $f(x)=x^4-\sqrt{m}$, 
hence in Lemma \ref{l1} $a_1=0,a_2=0,a_3=0,a_4=-\sqrt{m}$. We obtain
\begin{eqnarray}
F(U,V)&=&U(U^2+4\sqrt{m}V^2), \label{fuv}\\
Q_1(X,Y,Z)&=&X^2-\sqrt{m}Z^2, \label{q1}\\
Q_2(X,Y,Z)&=&Y^2-XZ.\label{q2}
\end{eqnarray}
Let $H=\{+1,-1\}$, in case $m=-1$ let $H=\{+1,-1,i,-i\}$.
Equation (\ref{res}) implies
\[
U(U^2+4\sqrt{m}V^2)=\varepsilon
\]
with $\varepsilon\in H$. This implies $U=\eta\in H$, hence
\[
V^2=\frac{\eta-\varepsilon^2}{4\sqrt{m}}.
\]
But $4\sqrt{m}$ can not divide $\eta-\varepsilon^2$, except if $\eta-\varepsilon^2=0$, whence $V=0$.
In view of (\ref{q12}) this implies 
\[
Q_1(X,Y,Z)=X^2-\sqrt{m}\, Z^2=\varepsilon,\;\;\; Q_2(X,Y,Z)=Y^2-XZ=0.
\]
If $Z=0$, then $Y=0$ and $X=\varepsilon_0\in H$ is a unit in $M$. This case provides possible solutions
\begin{equation}
X=\varepsilon_0,Y=Z=0.
\label{possible}
\end{equation}
If $Z\ne 0$, then we follow the arguments of \cite{g42} and parametrize $X,Y,Z$.
We have
\[
Q_0(X,Y,Z)=U\cdot Q_2(X,Y,Z)-V\cdot Q_1(X,Y,Z)=\varepsilon(Y^2-XZ)
\]
whence
\begin{equation}
Y^2-XZ=0.
\label{q0}
\end{equation}
A non-trivial solution of this equation is $X_0=1,Y_0=1,Z_0=1$.
Let $R,P,Q\in M, R\ne 0$ and write $X,Y,Z$ in the form
\begin{eqnarray}
X&=&X_0 R+P=R+P,\nonumber\\
Y&=&Y_0 R+Q=R+Q,\label{xyz}\\
Z&=&Z_0 R=R\nonumber.
\end{eqnarray}
Sustituting these representations of $X,Y,Z$ into (\ref{q0}) we obtain
\begin{equation}
(2Q-P)R=-Q^2.
\label{rr}
\end{equation}
We multiply all equations in (\ref{xyz}) by $2Q-P$ and replace $(2Q-P)R$ by $-Q^2$ in view of (\ref{rr}).
Then we obtain common multiples of $X,Y,Z$ in terms of quadratic forms of $P,Q$:
\begin{eqnarray}
SX&=&-P^2+2PQ-Q^2,\nonumber\\
SY&=&-PQ+Q^2,\label{pq}\\
SZ&=&-Q^2\nonumber
\end{eqnarray}
where $S\in M$. Following the arguments of \cite{g42} we multiply all equations in (\ref{pq}) by the square of a common
denominator of $P,Q$ and replace $S,P,Q$ by integer parameters $S,P,Q\in\Z_M$. Moreover $S\in\Z_M$ divides the determinant of the 
matrix
\[
\left(
\begin{array}{rrr}
-1&2&-1\\
0&-1&1\\
0&0&-1
\end{array}
\right)
\]
composed from the coefficients of $P^2,PQ,Q^2$ in the representation (\ref{pq}) of $X,Y,Z$.
Therefore we obtain that $S$ is a units in $M$.
Substituting the representations (\ref{pq}) of $X,Y,Z$ into $Q_1(X,Y,Z)=U$ we obtain
\begin{equation}
P^4-4P^3Q+6P^2Q^2-4PQ^3+(1-\sqrt{m})Q^4=S^2\cdot U=\rho,
\label{q1eq}
\end{equation}
where $\rho\in H$ is a unit. Finally, setting $P=A+B,Q=B$ we arrive to
\begin{equation}
A^4-\sqrt{m}\, B^4=\rho \;\;\; {\rm in}\;\; A,B\in\Z_M.
\label{ab}
\end{equation}

Note that by (\ref{100}) we have 
$\max(|X|,|Y|,|Z|)\le (1+\sqrt{5000})10^{200}=c_0$, therefore using (\ref{pq}) we have
\[
|Q^2|=|Z|\le c_0,\; |PQ|\le |Q^2|+|Y|\le 2c_0,\; |P|^2\le |X|+2|PQ|+|Q^2|\le 6c_0,
\]
hence $\max(|P|,|Q|)\le \sqrt{6c_0}$ and $\max(|A|,|B|)\le 2 \sqrt{6c_0}$,
finally
\[
\max(|a_1|,|b_1|,|a_2|,|b_2|)\le \max(|A|,|B|)\le 2 \sqrt{6c_0}\le 41.4856\cdot 10^{100}<10^{101}.
\]

\subsection{Calculating "small" solutions of the quartic relative Thue equation over $M$}

Let $A=a_1+\sqrt{m}\, b_1,B=a_2+\sqrt{m}\, b_2$ be a solution of (\ref{ab}). In the following
we determine all solutions of (\ref{ab}) with
\begin{equation}
C=\max(|a_1|,|b_1|,|a_2|,|b_2|)<10^{101}.
\label{kicsi}
\end{equation}
In order to make our algorithm efficient, in the following we sharpen the usual estamates
applying specialities of pure octic.

Denote by $\alpha^{(j)}\; (j=1,2,3,4)$ the relative conjugates of $\alpha=\sqrt[8]{m}$ over $M$.
Set $\beta^{(j)}=A-\alpha^{(j)}B\; (j=1,2,3,4)$.
Then equation (\ref{ab}) implies
\begin{equation}
\prod_{j=1}^4|\beta^{(j)}|\le 1.
\label{beta}
\end{equation}
Let $j_0$ be the index with $|\beta^{(j_0)}|=\min_{1\le j\le 4}|\beta^{(j)}|$, then 
\begin{equation}
|\beta^{(j_0)}|\le 1.
\label{j01}
\end{equation}
For $j\ne j_0$ we have
\begin{equation}
|\beta^{(j)}|\ge |\beta^{(j)}-\beta^{(j_0)}|-|\beta^{(j_0)}|\ge |\alpha^{(j)}-\alpha^{(j_0)}| \cdot |B|-1
\ge 0.9\cdot |\alpha^{(j)}-\alpha^{(j_0)}| \cdot |B|
\label{betaj}
\end{equation}
assumed 
\begin{equation}
|B|\ge\frac{10}{|\alpha^{(j)}-\alpha^{(j_0)}|}.
\label{B1}
\end{equation}
Therefore by (\ref{beta})  and  (\ref{betaj})  we obtain
\begin{equation}
|\beta^{(j_0)}|\le \frac{1}{\prod_{j\ne j_0}|\beta^{(j)}|}\le 
\frac{1}{0.9^3 \prod_{j\ne j_0}|\alpha^{(j)}-\alpha^{(j_0)}|}\cdot |B|^{-3}.
\label{j0}
\end{equation}
By (\ref{j01}) we obtain
\begin{equation}
|A|\le |\beta^{(j_0)}|+|\alpha^{(j_0)}|\cdot |B|\le 1.1 \cdot |\sqrt[8]{m}|\cdot |B|,
\label{aaa}
\end{equation}
assumed
\begin{equation}
|B|\ge\frac{10}{|\sqrt[8]{m}|}.
\label{B2}
\end{equation}
We have $\max(|a_1|,|b_1|)\le |A|$ and $\max(|a_2|,|b_2|)\le |B|$, hence by (\ref{aaa}) we obtain
\[
C=\max(|a_1|,|b_1|,|a_2|,|b_2|)\le 1.1 \cdot |\sqrt[8]{m}| \cdot |B|,
\]
whence
\begin{equation}
|B|^{-1}\le 1.1 \cdot |\sqrt[8]{m}| \cdot C^{-1},
\label{cc}
\end{equation}
which in view of (\ref{j0}) implies
\begin{equation}
|\beta^{(j_0)}|\le 
\frac{1}{0.9^3 \prod_{j\ne j_0}|\alpha^{(j)}-\alpha^{(j_0)}|}\cdot
\frac{(1.1\cdot |\sqrt[8]{m}|)^3}{C^3}
\le \left(\frac{1.1}{0.9}\right)^3\cdot\frac{1}{4}\cdot C^{-3}\le 0.4564\cdot C^{-3},
\label{jj00}
\end{equation}
using that $|\alpha^{(j)}-\alpha^{(j_0)}|$ is equal to $\sqrt{2}\cdot |\sqrt[8]{m}|$ for two values of $j\ne j_0$ and is equal to
 $2\cdot |\sqrt[8]{m}|$ for the third value.

\subsection{Reduction}

By (\ref{jj00}) we obtain
\begin{equation}
|(a_1+\sqrt{m}\, b_1)-\alpha^{(j_0)} (a_2+\sqrt{m}\, b_2)|\le d \cdot C^{-3}.
\label{redineq}
\end{equation}
We use this inequality with $d= 0.4564$ to reduce the upper bound (\ref{kicsi}) for $C$.
In view of (\ref{B1}), (\ref{B2}) this estimate is valid for $|B|\ge 10$. The values $a_2,b_2$
with $|a_2|,|b_2|\le 10$ will be tested separately.

Let $H$ be a large constant and 
consider the lattice generated by the columns of the matrix
\[
\left(
\begin{array}{cccc}
1&0&0&0\\
0&1&0&0\\
0&0&1&0\\
0&0&0&1\\
   H&H\Re(\sqrt{m})&H\Re(-\alpha^{(j_0)})&H\Re(-\alpha^{(j_0)}\sqrt{m})\\
   0&H\Im(\sqrt{m})&H\Im(-\alpha^{(j_0)})&H\Im(-\alpha^{(j_0)}\sqrt{m})\\
\end{array}
\right).
\]

\vspace{0.5cm}

\begin{lemma} (Lemma 5.3 of \cite{book})
Denote by $\ell_1$ the first vector of the LLL reduced basis of this lattice.
If $C\le C_0$ and $H$ is large enough to have
\[
|\ell_1|\ge \sqrt{40}\cdot C_0,
\]
then
\[
C\le \left(\frac{d\cdot H}{C_0}\right)^{1/3}.
\]
\label{redlemma}
\end{lemma}

For a certain $C_0$ the suitable $H$ is of magnitude $C_0^2$, therefore 
the reduced bound is about $C_0^{1/3}$. To make our procedure stable we were mostly 
using $H=10\cdot C_0^2$. A typical
sequence of reduced bounds starting from $C_0=10^{101}$ was the following:
\[
\begin{array}{|c|c|c|c|c|c|}
\hline
C            &10^{101}            &7.6992\cdot 10^{33}    &3.2754\cdot 10^{11}  &11434              &37 \\ \hline
H            &10^{203}            &5.9278\cdot 10^{68}    &1.0728\cdot 10^{24}  &1.3073\cdot 10^{9} &8556.25\\   \hline
{\rm new}\; C&7.6992\cdot 10^{33} &3.2754\cdot 10^{11}    &11434                &37                 &4\\ \hline
\end{array}
\]

Note that this reduction procedure must be performed for $j_0=1,2,3,4$.
for each considered $m$. The reduced bound was under 10, except for 
$m=-5,-62,-197$ when we obtained the reduced bounds $42,18,13$, respectively.

\subsection{Testing small solutions}

For each considered $m$ we tested if for $0\le b_1\le 10,-10\le b_2\le 10$ there exists $a_1,b_1$
satisfying (\ref{ab}). (For $m=-5,-62,-197$ we used the bounds $42,18,13$, respectively.)
Calculating $B=a_2+\sqrt{m}\, b_2$ and $G=\sqrt{m}B^4+\rho$ (for all units $\rho$ in $M$) we
checked if any of the fourth roots of $G$ are of type $a_1+\sqrt{m}\, b_1$ with integer values
$a_1,b_1$. This can be easily done by separating the real and complex parts of $\sqrt[4]{G}$.
If $a_2=b_2=0$ then we obtained the solution $a_1=\pm 1,b_1=0$, for $m=-1$ in addition 
$a_1=0,b_1=\pm 1$. If $a_2^2+b_2^2>0$ then we found no solutions except
$a_1=0,b_1=0,a_2=0,b_2=\pm 1$ and $a_1=0,b_1=0,a_2=\pm 1,b_2=0$ for $m=-1$.
In view of $P=A+B,Q=B$ and (\ref{pq}) these tuples gave only solutions with $Z=0$, except for
$X=0,Y=0,Z=\pm 1$ for $m=-1$.

If $Z=0$ then in view of (\ref{possible}) we obtain $X=\pm 1,Y=0,Z=0$ and in case $m=-1$
in addition $X=\pm i,Y=0,Z=0$.

Having the possible values of $X,Y,Z$, we know that up to $\Z_M$-equivalence  all
possible generators of relative power integral bases of $K$ over $M$ are of the form
\[
\delta=X\alpha+Y\alpha^2+Z\alpha^3.
\]
This implies that the possible generators of (absolute) power integral bases of $K$ 
(over $\Q$) are up to $\Z$-equivalence of the form
\[
\gamma=y_0\sqrt{m}+\varepsilon \delta
\]
where
$\varepsilon$ is a unit in $M$. For each $m$ and for all possible $\delta,\varepsilon$ 
the equation $I(\gamma)=1$ gives a polynomial equation of degree 16 in $y_0$
with integer coefficients. Solving these equations we can easily check if for 
the given $\delta$ there exists a corresponding $y_0$.
Summarizing, we obtained that up to $\Z$-equivalence the only generator of power integral bases 
of $K$ with (\ref{100}) is $\alpha$, except for $m=-1$, when we additionally have 
$i\alpha$, $\alpha^3$, $i\alpha^3$.

\section{Computational remarks}

All calculations were performed in Maple \cite{maple}.
The reduction procedure took 100 minutes, testing small values  of $a_2,b_2$
took 60 minutes, searching for suitable $y_0$ values required 25 minutes,
all together 185 minutes for the 2024 fields.

\end{document}